\documentclass{article}
\usepackage{emsnewsletter}

\usepackage[english]{babel}
\usepackage{hyperref}       
\usepackage{url}            
\usepackage{booktabs}       
\usepackage{amsfonts}       
\usepackage{tablefootnote}
\usepackage{verbatim}
\usepackage{nicefrac}       
\usepackage{microtype}      
\usepackage{lipsum}
\usepackage{mathtools}

\usepackage{amsmath,amssymb,amsthm}
\usepackage{algorithm,algorithmic}
\usepackage{pifont}
\usepackage{cases}
\usepackage{subcaption,graphicx}
\usepackage{stackengine}    
\usepackage{wrapfig}
\usepackage{enumitem}
\usepackage{multirow}
\usepackage{xcolor}
\usepackage{array}

\renewcommand{\le}{\leqslant}
\renewcommand{\ge}{\geqslant}

\DeclareMathOperator{\image}{Im}

\newcommand{\R}{\mathbb{R}}

\newcommand{\E}{\mathbb{E}}

\def\<#1,#2>{\langle #1,#2\rangle}

\newcommand{\avg}[1]{{\color{black}#1}}

\newcommand{\pd}[1]{{\color{black}#1}}

\newcommand{\dd}[1]{{\color{black}#1}}

\newcommand{\eg}[1]{{\color{black}#1}}

\newcommand{\gea}[1]{{\color{black}#1}}

\newcommand{\al}[1]{{\color{black}#1}}

\author{
	\textit{Alexander Gasnikov$^{1,2,3}$,
	Darina Dvinskikh$^{1,4,6}$,  Pavel Dvurechensky$^7$, 
	Eduard Gorbunov$^5$,  Aleksandr Beznosikov$^{1}$, Aleksandr Lobanov$^{1,8}$}\\
	~1. Moscow Institute of Physics and Technology, Moscow, Russia \\
	~2. Caucasus Mathematical Center, Adyghe State University, Maikop, Russia\\
 	~3. \al{Innopolis University, Innopolis, Russia}\\
	~4. IITP RAS, Moscow, Russia\\
	~5. MBZUAI, Abu Dhabi, UAE\\
	~6. HSE University, Moscow, Russia\\
	~7. WIAS, Berlin, Germany\\
 	~8. \al{ISP RAS Research Center for Trusted Artificial Intelligence, Moscow, Russia} 
}
\title{\Large{Randomized gradient-free methods in convex optimization}}

\begin{document}
   
\maketitle



\pd{
}\\


\section{Introduction}
\gea{C}onsider convex optimization problem
\begin{align}\label{eq:problem}
    \min_{x\in Q \subseteq \R^d}~ f(x)
\end{align}
\pd{with convex feasible set $Q$, convex objective $f$ possessing the} zeroth-order (gradient/derivative-free) oracle \cite{rosenbrock1960automatic}. \pd{The latter} means that \gea{one has} an access only to \pd{the values of the objective $f(x)$} rather than \pd{to its gradient} $\nabla f(x)$ that is more popular for numerical methods \cite{polyak1987introduction,nesterov2018lectures}. 
\gea{One can find motivating examples in \cite{granichin2003randomizirovannye,spall2005introduction,conn2009introduction,larson2019derivative}.}
\gea{This paper focuses} on numerical zeroth-order methods based on first-order methods with different approximation models of the gradient.

The first such methods (\pd{often referred to as Kiefer--Wolfowitz methods}) were proposed about 70 years ago (see, e.g., \cite{wasan2004stochastic} \pd{for historical remarks})  and further developed in the \pd{following} 20--30 years \cite{ermoliev1976stochastic,nemirovsky1983problem,polyak1987introduction,spall2005introduction}.

\dd{The} interest in \dd{gradient-free} methods has grown considerably \pd{since mid-2000s} \cite{agarwal2010optimal,agarwal2011stochastic,duchi2012randomized,duchi2015optimal,gasnikov2016gradient,gasnikov2017stochastic,nesterov2017random,shamir2017optimal,gorbunov2022accelerated}. Particularly, in modern papers \cite{gasnikov22a,dvinskikh2022gradient,akhavan2022gradient}, \pd{the} authors try to develop the best known gradient-free algorithms at least for one of three main criteria: oracle complexity \pd{(number of zeroth-oracle calls to guarantee certain accuracy)}, iteration complexity \pd{(number of iterations to guarantee certain accuracy)}, maximum level of admissible noise \pd{still allowing to guarantee certain accuracy}. \gea{This survey focuses on the same aspects and presents the algorithms through the prism of optimality according to the listed criteria.}

\section{Full gradient approximation. Smooth case}\label{sec:direct}
\gea{Assume that one can query zeroth-order oracle which returns }
\pd{a noisy} value of $f(x)$: $f_{\delta}(x) = f(x) +\delta(x)$, where  \dd{the noise is bounded by some $\Delta\geqslant 0$, i.e.,} $\left|\delta(x)\right| \le \Delta$. In some applications, (see e.g. \cite{bogolubsky2016learning}) the \pd{larger} is $\Delta$ the cheaper is the call of inexact oracle $f_{\delta}$. If $f$ in \eqref{eq:problem} has $L$-Lipschitz gradient (e.g., \pd{w.r.t. the Euclidean} $2$-norm)\pd{, then \gea{one can} approximate partial derivatives using the forward differences of inexact values of the objective} as follows:
$$\left|\frac{\partial f(x)}{\partial x_i} -\frac{f_{\delta}(x + h e_i) - f_{\delta}(x)}{h}\right| \le \frac{Lh}{2}+\frac{2\Delta}{h}\le 2\sqrt{L\Delta}\pd{,}
$$
where \pd{$e_i$ is the $i$-th basis vector} 
and $h = 2\sqrt{\Delta/L}$.

\dd{Furthermore,} if $f$ in \eqref{eq:problem} has $\bar{L}$-Lipschitz Hessian (e.g., \pd{w.r.t. the} $2$-norm) then \gea{one can} \pd{approximate partial derivatives using symmetric differences} as follows:
$$\left|\frac{\partial f(x)}{\partial x_i} -\frac{f_{\delta}(x + h e_i) - f_{\delta}(x - he_i)}{2h}\right| \le \frac{\bar{L}h^2}{6}+\frac{\Delta}{h}\le 2\bar{L}^{1/3}\Delta^{2/3}\pd{,}
$$
where $h = \left(3\Delta/\bar{L}\right)^{1/3}$.
Similarly, if $f$ has Lipschitz \pd{$n$-th} order derivatives, one can approximate partial derivatives by using \pd{$n$-th} order finite-difference scheme with accuracy $\sim \Delta^{\frac{n}{n+1}}$ \cite{vasin2021stopping,berahas2022theoretical}. \pd{This} means that \gea{one can} build a finite-difference inexact gradient model $\tilde{\nabla} f(x)$ with accuracy   
\begin{equation}\label{eq:err}
\|\nabla f(x) -  \tilde{\nabla} f(x)\|_2 \lesssim \sqrt{d}\Delta^{\frac{n}{n+1}}.
\end{equation}
\pd{Here and below $\lesssim$ is used to denote inequality up to possible dimension-dependent constants.}

Most of first-order\pd{, i.e., gradient,} algorithms under a proper stopping rule are <<robust>>, i.e., \dd{they} do not accumulate additive error in the gradient \cite{juditsky2011first,gorbunov2019optimal,vasin2021stopping}. \pd{Particularly}, if \pd{the goal is}  to solve \eqref{eq:problem} with a sufficiently small accuracy $\varepsilon$ in \pd{terms of the objective} function value (\pd{i.e., s.t. the point $\bar{x}^N$ generated after $N$ steps satisfies} $f\left(\bar{x}^N\right) - f(x_*) \le \varepsilon$), \pd{then} the additive error should be proportional to $\varepsilon$, i.e., $\sqrt{d}\Delta^{\frac{n}{n+1}} \simeq \varepsilon/R,$ where $R = \|x^0 - x_*\|_2$ with $x^0$ being the starting point and  $x_*$ being the solution to \eqref{eq:problem} closest to $x^0$. Note that it is easy to show that $R$ can be bounded by the diameter of $Q$. \pd{Showing} that $R$ is determined as a distance between the starting point $x^0$ and the solution $x_*$ \pd{requires} a more accurate analysis. This remark is important when $Q$ is not a compact set or if \pd{the goal is} to generalize the results for \pd{minimizing} strongly convex \pd{functions} by using the restart-technique \cite{dvinskikh2022gradient}. \pd{Thus,} the maximum admissible level of noise is $\Delta \lesssim \varepsilon^{\frac{n+1}{n}}d^{-\frac{n+1}{2n}}$. If $n \to \infty$, this bound has a form $\Delta \lesssim \frac{\varepsilon}{\sqrt{d}}$, \pd{which is} a lower bound \cite{singer2015information,risteski2016algorithms}. \dd{However}, \pd{as it is explained} below (see Section~\ref{sec:non_smooth_l2}), \pd{for finite values of $n$, the bound for the} admissible level of noise can be improved.

\pd{In the case when $f$ is an analytical function and $\Delta = 0$,} \gea{one can} use \pd{a complex} representation of the partial derivatives \avg{\cite{squire1998using,jongeneel2021small}} 
$$
\left|\frac{\partial f(x)}{\partial x_i} -\image\frac{f(x+ih)}{h} \right| \lesssim h^2.
$$
Typically in practice, \pd{when the floating point arithmetic is used,} \gea{one has} a double precision $\tilde{\Delta} \simeq 1.1\cdot 10^{-16}$ (so $h$ should be smaller than $10^{-8}$), and with this precision, \pd{when calculating} $\partial f(x)/\partial x_i$, \gea{one has} $\sqrt{d}\tilde{\Delta}$ error in the gradient, \pd{meaning that} \gea{one can} solve \eqref{eq:problem} with accuracy $\varepsilon$ proportional to this error. A close result for large $n$ can be obtained by using the finite-difference approach with \pd{the change} $\tilde{\Delta}\to\Delta$. This example \gea{is described} to \pd{emphasise} the difference between the inexactness \pd{of the zeroth-order oracle} and the finite-precision inexactness. The last one is less restrictive and opens up opportunities to use more unconventional approaches for gradient approximation.

\gea{Now, everything it ready for a description of} the most simple way to build optimal gradient-free methods. The basic observation here is as follows: optimal Nesterov's accelerated gradient method \cite{nesterov2018lectures} is robust \cite{gorbunov2019optimal,vasin2021stopping}. Starting from this gradient-based algorithm \gea{one can} build an optimal gradient-free algorithm \cite{nemirovsky1983problem} by replacing the \pd{exact} gradient by its finite-difference \pd{approximation}. \pd{Such} gradient-free method, in general, is optimal \pd{in terms of the number of} oracle calls $O\left(d\sqrt{\frac{LR^2}{\varepsilon}}\right)$ and iteration complexity $O\left(\sqrt{\frac{LR^2}{\varepsilon}}\right)$, and close to optimal algorithms (but is sill not optimal) in terms of the maximum level of admissible noise \pd{with the bound} $\Delta \lesssim \varepsilon^{\frac{n+1}{n}}d^{-\frac{n+1}{2n}}$. \pd{The lower bound is not known}, but the best \pd{known result}  $\Delta \lesssim \varepsilon^{\frac{n+2}{n+1}}d^{-\frac{1}{2}}$
requires more iterations or oracle calls, 
see Section~\ref{sec:non_smooth_l2} for $n=1$.
 
\section{Randomized gradient approximation. Smooth case}\label{sec:rand_smooth}
At the end of Section~\ref{sec:direct}, a rather simple approach that \pd{is} based \pd{on the} gradient approximation \gea{is described}. This approach was announced to be optimal in terms of the number of oracle calls. In this section, alternative randomized approximations that can improve oracle complexity under some additional assumptions \gea{are considered}.

\subsection{Coordinate-wise randomization}\label{sec:coord}
Coordinate methods have a long history go\dd{ing} back to \dd{the} Gauss--Seidel method \cite{polyak1987introduction}. A contemporary view on these methods was proposed by Yu.~Nesterov \cite{nesterov2012efficiency}, \pd{followed by a large number of works making this area of optimization well-developed} \cite{richtarik2014iteration,wright2015coordinate,gasnikov2015accelerated,nesterov2017efficiency,hanzely2020variance}. The core idea of these methods is to replace \pd{exact} gradient by \pd{its unbiased approximation} 
$$
\tilde{\nabla} f(x,i) = d\cdot\frac{\partial f(x)}{\partial x_i} \cdot e_i,
$$
where $i$ is chosen from $1,...,d$ \pd{randomly} with equal probabilities. If the oracle returns \pd{the exact value of} $\partial f(x)/\partial x_i$, then $\mathbb{E}_{i} \tilde{\nabla} f(x,i) = \nabla f(x)$ and  $\mathbb{E}_{i} \|\tilde{\nabla} f(x,i)\|_2^2 = d\|\nabla f(x)\|_2^2$. \pd{This leads to the} <<strong growth condition>> \cite{vaswani2019fast} and the same oracle calls complexity bound \pd{as that of for the methods with full gradient approximation discussed at the end of Section~\ref{sec:direct}, i.e., } $O\left(d\sqrt{\frac{LR^2}{\varepsilon}}\right)$. \pd{Here and below, when considering randomized algorithms, $\varepsilon$ is a desired accuracy in terms of expectation of the objective function  value.} A more accurate analysis allows to improve this bound. Namely, let $L_i$ \pd{satisfy} 
$$
\left| \frac{\partial f(x + he_i)}{\partial x_i} - \frac{\partial f(x)}{\partial x_i} \right| \le L_i h
$$
and $\bar{L} = \frac{1}{d}\sum_{i=1}^d L_i$. For example, if $f(x) = \frac{1}{2}\langle Ax,x \rangle - \langle b,x\rangle$,  $L = \lambda_{\max}(A)$, $\bar{L} = \frac{1}{d}\text{tr}(A) = \frac{1}{d}\sum_{i=1}^d \lambda_i$. Hence, $\frac{1}{d}L\le \bar{L}\le L$  \cite{nesterov2017efficiency}.
A proper version of accelerated (catalyst-based) Nesterov's coordinate descent method \cite{ivanova2021adaptive,ivanova2022oracle} requires  $\tilde{O}\left(d\sqrt{\frac{\bar{L}R^2}{\varepsilon}}\right)$ oracle calls that could be better in $\tilde{O}\left(\sqrt{d}\right)$ times than the previous bound. \gea{Further discussions of the gradient-free modifications of coordinate descents are given in} Section~\ref{sec:ZOsmooth}.

\subsection{Random search randomization}\label{sec:l2}
The core idea of random-search methods (see, e.g., \cite{ermoliev1976stochastic,dvurechensky2021accelerated}) is to replace \dd{the} \pd{exact} gradient by \pd{its unbiased random approximation}
\begin{equation}\label{eq:l2_rand}
\tilde{\nabla} f(x,e) = d\cdot\langle \nabla f(x),e\rangle \cdot e,
\end{equation}
where $e$ is chosen \pd{uniformly} at random from the unit \dd{E}uclidean sphere $S^d_2(1)$. If the oracle return\pd{s} \pd{the exact value of}  $\langle \nabla f(x),e\rangle$, then $\mathbb{E}_{e} \tilde{\nabla} f(x,e) = \nabla f(x)$ and  $\mathbb{E}_{e} \|\tilde{\nabla} f(x,e)\|_2^2 = d\|\nabla f(x)\|_2^2$. \dd{Thus,} the <<strong growth condition>> \cite{vaswani2019fast} \gea{still holds} and \gea{one} can \pd{obtain} the same oracle calls complexity \pd{as that of for the methods with full gradient approximation discussed at the end of Section~\ref{sec:direct}, i.e., } $O\left(d\sqrt{\frac{L\|x^0-x_*\|_2^2}{\varepsilon}}\right)$. A more accurate analysis of proper accelerated (linear-coupling) random-search method \cite{dvurechensky2021accelerated} allows \dd{us} to improve this bound \pd{to} $\tilde{O}\left(\sqrt{d\frac{L\|x^0-x_*\|_1^2}{\varepsilon}}\right)$, where $L$ (as it was above) is a Lipschitz gradient constant of $f$ in \dd{the} $2$-norm. \pd{On the one hand,} $\|x^0-x_*\|_2^2\le\|x^0-x_*\|_1^2 \le d\|x^0-x_*\|_2^2$ and if $\|x^0-x_*\|_1^2\simeq d \|x^0-x_*\|_2^2$, \pd{there is no improvement in the new bound}. \pd{On the other hand,} if $\|x^0-x_*\|_2^2 \simeq \|x^0-x_*\|_1^2$, the oracle calls complexity \pd{improves by the factor of} $\tilde{O}\left(\sqrt{d}\right)$.

\subsection{Gradient-free randomized methods. Smooth case}\label{sec:ZOsmooth}
In Sections~\ref{sec:coord} and \ref{sec:l2}, the ways that \dd{make it possible} to improve oracle complexity (\pd{but worsen the} iteration complexity at the same time) \gea{are described}. \pd{But, these approaches make use of the} partial derivatives oracle. \pd{These results can be naturally generalized to the setting of the zeroth-oracle inexact oracle by the following modification of example \eqref{eq:l2_rand}:}
\begin{equation}\label{eq:l2_rand_GF}
\tilde{\nabla} f(x,e) = d \frac{f_{\delta}(x+he) - f_{\delta}(x-he)}{2h} e.
\end{equation}
If $f$ has Lipschitz $n$-th order derivatives, the <<kernel-based>> approximation from \cite{polyak1990optimal,bach2016highly,bubeck2017kernel,akhavan2020exploiting,novitskii2022improved,lobanov2023highly,lobanov2023black} \pd{can be used} to avoid $n$-th order finite-difference schemes:
\begin{equation}\label{eq:l2_rand_GF}
\tilde{\nabla} f(x,e) = d \frac{f_{\delta}(x+h\kappa e) - f_{\delta}(x-h\kappa e)}{2h}K_p(\kappa)e ,
\end{equation}
where $\kappa$ is uniformly distributed in $[-1,1]$, $e$ is uniformly distributed in $S^d_2(1)$, $\kappa$ and $e$ are independent, $K_p(\kappa)$ is a kernel function that satisfies $\E[K_n(\kappa)] = 0, \,
    \E[\kappa K_n(\kappa)] = 1, \, 
    \E [\kappa^j K_n(\kappa)] = 0,\, j=2, \dots, n, \,
    \E\left[ |\kappa|^{n+1}|K_n(\kappa)|\right] < \infty$.
For example, 
\begin{center}
$K_n(\kappa) = 3\kappa$, $n = 1,2$; $K_n(\kappa) = \frac{15\kappa}{4}\left(5 - 7\kappa^2\right)$, $n = 3, 4$ ... \end{center}

A more accurate analysis allows \pd{one} to replace the gradient approximation condition \eqref{eq:err} by \pd{the following} \cite{akhavan2020exploiting,gasnikov22a,novitskii2022improved}: for all $r \in  \mathbb{R}^d$,
\begin{equation}\label{eq:unbias}
\mathbb{E}_e\left\langle \tilde{\nabla} f(x, e)  - \nabla f(x), r\right\rangle \lesssim \sqrt{d}\left(h^n+\frac{\Delta}{h}\right)\|r\|_2 \le \sqrt{d}\Delta^{\frac{n}{n+1}}\|r\|_2.
 \end{equation}
For coordinate-wise randomization, \gea{one can} further improve this bound  replacing $\sqrt{d}\|r\|_2$ \dd{by} $\|r\|_1$. If \gea{one} wants to solve \eqref{eq:problem} with a sufficiently small accuracy $\varepsilon$ in \pd{terms of} \dd{the} function value, i.e. $\E f\left(\bar{x}^N\right) - f(x_*) \le \varepsilon$ (this criteria \gea{will be called}  <<$\varepsilon$-suboptimality in \dd{the} expectation>> \gea{in this survey}), th\dd{e}n \dd{the} RHS  of \eqref{eq:unbias} should be of \pd{the} order \pd{of} $\varepsilon: \sqrt{d}\Delta^{\frac{n}{n+1}}R \simeq \varepsilon,$ where $R = \|x^0 - x_*\|_2$ \cite{gasnikov22a}. \pd{Thus,} the maximum admissible level of noise is $\Delta \lesssim \varepsilon^{\frac{n+1}{n}}d^{-\frac{n+1}{2n}}$, \pd{i.e., the same as it was for the algorithms using the full gradient} approximation, see Section~\ref{sec:direct}.

\section{Randomized gradient approximation. Non-smooth case}\label{sec:non_smooth}
If $f$ in \eqref{eq:problem} is no\dd{n-}smooth but is Lipschitz continuous in $p\in[1,2]$ norm:
 $|f(y) - f(x) | \leq M \|y - x \|_p$ (if $p=2$ this constant \gea{will be denoted} as $M_2$), then the finite-difference approximation approaches do not work. To \pd{resolve} this problem, different smoothing techniques are used \cite{nemirovsky1983problem,spall2005introduction,gasnikov2016gradient,nesterov2017random,akhavan2022gradient,lobanov2022,lobanov2023zero,lobanov2023accelerated,Lobanov_OPTIMA23,kornilov2023accelerated}. The Gaussian smoothing \cite{nesterov2017random} \gea{is not considered} because it does not give anything new compared to smoothing on the $l_2$-ball \cite{vasin2021stopping}. \gea{This} section \gea{is considered} to be the most important part of the survey. 
 

\subsection{$l_2$-ball smoothing scheme}\label{sec:non_smooth_l2}
This technique goes back to the works of A.~Nemirovski, Yu.~Ermoliev, A.~Gupal of mid- and late \pd{1970s} \cite{nemirovsky1983problem,granichin2003randomizirovannye}, see also \cite{spall2005introduction}. \pd{In this section,} \gea{the narrative} follows  \cite{gasnikov22a}.

The first element of the approach \cite{gasnikov22a} is the randomized \textit{smoothing} for the non-smooth objective $f$:
$$
f_{\gamma}(x)=  \mathbb{E}_{u} f(x + \gamma u),
$$
where $u \sim RB^d_2(1)$, i.e., $u$ is random vector uniformly distributed on a unit euclidean ball $B^d_2(1)$. The following theorem is a generalization of the results \cite{YOUSEFIAN201256,duchi2012randomized} for non-Euclidean norms.

\begin{theorem}[properties of $f_{\gamma}$]
\label{main_properties} For all $x,y\in Q$
\begin{itemize}[noitemsep,topsep=0pt,leftmargin=10pt]
    \item $f(x) \leq f_{\gamma}(x) \leq f(x) + \gamma M_{2}$;
    \item $f_{\gamma}(x)$ has $L = 
    \dfrac{\sqrt{d} M}{\gamma} $-Lipschitz gradient:
    $$\|\nabla f_{\gamma}(y) - \nabla f_{\gamma}(x) \|_q \leq L \| y - x\|_p,$$ 
     where $q$ is such that $1/p + 1/q = 1$. 
\end{itemize}
\end{theorem}

The second very important element of the approach \cite{gasnikov22a} goes back to O.~Shamir \cite{shamir2017optimal}, who proposed a special unbiased stochastic gradient for $f_{\gamma}$ with small variance (note that the forward finite-difference has large variance \cite{duchi2015optimal}):
\begin{equation}\label{sg}
    \nabla f_{\gamma}(x,e) = d\frac{f(x+\gamma e) - f(x-\gamma e)}{2\gamma}e,
\end{equation}
where $e \sim RS^d_2(1)$ is a random vector uniformly distributed on $S^d_2(1)$ (see also \eqref{eq:l2_rand_GF} for comparison). To simplify the further derivations, here and below (until the end of the section) \gea{it is assumed} that $\Delta = 0$.

The following theorem is a combination of the results from \cite{shamir2017optimal,gorbunov2019upper,beznosikov2020derivative}.
    \begin{theorem}[properties of $\nabla f_{\gamma}(x,e)$]
\label{main_properties2} 
For all $x\in Q$
 \begin{itemize}[noitemsep,topsep=0pt,leftmargin=10pt] 
    \item $\nabla f_{\gamma}(x,e)$ is an unbiased approximation for $\nabla f_{\gamma}(x)$:
    $$\mathbb{E}_e \left[\nabla f_{\gamma}(x, e) \right] = \nabla f_{\gamma}(x);$$
    
    \item  $\nabla f_{\gamma}(x,e)$ has bounded variance (second moment):
   \begin{equation}\label{eq:var}
    \mathbb{E}_e \left[ \| \nabla f_{\gamma}(x, e) \|^2_q \right] \le \sqrt{2} \min \left\{ q, \ln d \right\} d^{\frac{2}{q}} M_2^2,
    \end{equation}
   where $1/p + 1/q = 1$. 
   \end{itemize}
   \end{theorem}
   
Based on the two elements above, \gea{everything is ready for the description of} a general approach, which for shortness, \gea{is called in this paper} as the \textit{smoothing scheme} (technique).

Assume that \gea{one has an access to} some batched algorithm \textbf{A}($L,\sigma^2$) that solves problem \eqref{eq:problem} under the assumption that $f$ is smooth and satisfies
\begin{equation*}\label{LL}
\|\nabla f(y) - \nabla f(x) \|_q \leq L \| y - x\|_p, 
\end{equation*}
and by using a stochastic first-order oracle that depends on a random variable $\eta$ and returns at a point $x$ an unbiased stochastic gradient $\nabla_x f(x,\eta)$ with bounded variance:  
\begin{equation*}
\label{sigma^2}
     \mathbb{E}_{\eta} \left[ \| \nabla_x f(x, \eta) - \nabla f(x)\|^2_q \right]\le \sigma^2.
\end{equation*}
Further, assume that, to reach $\varepsilon$-suboptimality in expectation, this algorithm requires $N(L,\varepsilon)$ successive iterations and $T(L,\sigma^2,\varepsilon)$ stochastic first-order oracle calls, i.e., \textbf{A}($L,\sigma^2$) allows batch-parallelization with the average batch size $$B(L,\sigma^2,\varepsilon) = \frac{T(L,\sigma^2,\varepsilon)}{N(L,\varepsilon)}.$$

\pd{The approach of \cite{gasnikov22a} is to apply} \textbf{A}($L,\sigma^2$) to the smoothed problem
 \begin{equation}\label{sm_problem}
    \min_{x\in Q\subseteq \mathbb{R}^d} f_{\gamma}(x)
\end{equation}
with 
\begin{equation}\label{gamma}
  \gamma = \varepsilon/(2M_{2})  
\end{equation}
and $\eta = e$, $\nabla_x f(x, \eta) = \nabla f_{\gamma}(x, e)$, where $\varepsilon > 0$ is the desired accuracy for solving problem \eqref{eq:problem} in terms of the suboptimality expectation.

According to Theorem~\ref{main_properties}, an $(\varepsilon/2)$-solution to \eqref{sm_problem} is an $\varepsilon$-solution to the initial problem \eqref{eq:problem}.
According to Theorem~\ref{main_properties} and \eqref{gamma}, \gea{the following inequality holds }
\begin{equation}\label{L}
   L \le \frac{2\sqrt{d}M M_{2}}{\varepsilon}, 
\end{equation}
and \gea{Theorem~\ref{main_properties2} implies} 
\begin{equation}\label{var}
  \sigma^2 \le \sigma^2(p,d) = 2\sqrt{2} \min \left\{ q, \ln d \right\} d^{\frac{2}{q}} M_2^2.  
\end{equation}

Thus, \gea{it is shown} that \textbf{A}($L,\sigma^2$) implemented using stochastic gradient \eqref{sg} is a zeroth-order method for solving non-smooth problem \eqref{eq:problem}. Moreover, to solve problem \eqref{eq:problem} with accuracy $\varepsilon>0$ this method suffices to make 
$$
N\left(\frac{2\sqrt{d}M M_{2}}{\varepsilon},\varepsilon\right) \; \text{ successive iterations and }
$$
$$
2T\left(\frac{2\sqrt{d}M M_{2}}{\varepsilon},\sigma^2(p,d),\varepsilon\right) \; \text{ zeroth-order oracle calls. }
$$

\gea{It is important to mention} that this approach is flexible and generic as \gea{one can} take different algorithms as \textbf{A}($L,\sigma^2$). For example, if \gea{one takes} batched Accelerated gradient method \cite{cotter2011better,lan2012optimal,devolder2013exactness,dvurechensky2016stochastic,gorbunov2019optimal}, then from \eqref{L}, \eqref{var} \gea{it follows} that
\begin{center}
$N(L,\varepsilon) = O\left(\sqrt{\frac{LR^2}{\varepsilon}}\right)=O\left(\frac{d^{1/4}\sqrt{M M_2}R}{\varepsilon}\right),$ 
$T(L,\sigma^2,\varepsilon) = O\left(\max\left\{N(L,\varepsilon),\frac{\sigma^2(p,d)R^2}{\varepsilon^2}\right\}\right)$
$=O\left(\frac{d^{2/q} M_2^2R^2}{\varepsilon^2}\right),$ 
\end{center}
where  $R = O\left(\|x^0 - x_*\|_p\sqrt{\min\left\{q,\ln d\right\}}\right)$. The last equality on $T$ assumes also that $\varepsilon \lesssim d^{-1/4}M_2^{3/2}R/M^{1/2}$ when $p=1$.

\begin{theorem}\label{FGM}
Based on the batched Accelerated gradient method, the Smoothing scheme applied to non-smooth problem \eqref{eq:problem}, provides a gradient-free method with $$
O\left(\frac{d^{1/4}\sqrt{M M_2}R}{\varepsilon}\right) \; \text{ successive iterations and}
$$ 
$$O\left(\frac{\min \left\{ q, \ln d \right\} d^{\frac{2}{q}} M_2^2 R^2}{\varepsilon^2}\right)=
\begin{cases}
   O\left(\frac{d M_2^2R^2}{\varepsilon^2}\right), \;  p = 2\\
   O\left(\frac{(\ln d) M_2^2R^2}{\varepsilon^2}\right), \; p = 1.
 \end{cases}
 $$
zeroth-order oracle calls, 
where $1/p + 1/q = 1$.
\end{theorem}

The simplest and the most illustrative result obtained in Theorem~\ref{FGM} from \cite{gasnikov22a} is as follows. Theorem~\ref{FGM} assumes that \gea{one can} solve problem \eqref{eq:problem} in the Euclidean geometry in $O\left(d^{1/4}M_2R/\varepsilon\right)$ successive iterations, each requiring $O\left(d^{3/4}M_2R/\varepsilon\right)$ oracle calls per iteration that can be made in parallel to make the total working time smaller. The dependence $\sim d^{1/4}/\varepsilon$ corresponds to the first part of the lower bound for iteration complexity \cite{diakonikolas2020lower,bubeck2019complexity}  $\min\left\{d^{1/4}\varepsilon^{-1}, d^{1/3}\varepsilon^{-2/3}\right\}$ (this bound assumes that $d$ is large enough, otherwise centre of gravity type estimate $d\ln \varepsilon^{-1}$ is better, see also the end of the Section~\ref{sec:non_smooth_one_point}). Note that  this lower bound is obtained for a wider class of algorithms  that allow $\text{Poly}\left(d,\varepsilon^{-1}\right)$ oracle calls per iteration. In particular, the authors of \cite{bubeck2019complexity} propose a close technique to the described above with an accelerated higher-order method (see, e.g., survey \cite{kamzolov2022exploiting} and references therein) playing the role of \textbf{A}($L,\sigma^2$). For the particular setting of $p=2$, they obtain a better bound $N\sim d^{1/3}/\varepsilon^{2/3}$ for the number of iterations in some regimes. Yet, they have significantly worse oracle complexity $T$, which makes it unclear how to use their results in practice. On the contrary, the described above approach from \cite{gasnikov22a} makes $T$ optimal simultaneously \cite{nemirovsky1983problem,duchi2015optimal}. To conclude, the described above approach from \cite{gasnikov22a} is optimal for iteration complexity and oracle complexity. But what is about maximum admissible level of noise? 

To answer this question, \gea{one} should note that in case of $\Delta > 0$ the biased inequality \eqref{eq:unbias} becomes
\begin{equation}\label{eq:unbias_non_smooth}
\mathbb{E}_e\left\langle \tilde{\nabla}  f_{\gamma}(x,e) - \nabla f(x), r\right\rangle \le \frac{\sqrt{d}\Delta}{\gamma}\|r\|_2,
 \end{equation}
and \eqref{eq:var} should be replaced by
\begin{equation}\label{eq:var_noise}
\mathbb{E}_e \left[ \| \nabla f_{\gamma}(x, e) \|^2_q \right] \lesssim \min \left\{ q, \ln d \right\} d^{\frac{2}{q}} \left(M_2^2 + \frac{d\Delta^2}{\gamma^2}\right).
\end{equation}
Since $\gamma = \varepsilon/(2M_{2}) $ (see \eqref{gamma}), from \eqref{eq:unbias_non_smooth} the following condition for $\Delta$ \gea{holds} (RHS of \eqref{eq:unbias_non_smooth} is compared with $\varepsilon$):
\begin{equation}\label{lon}
\Delta\lesssim \frac{\varepsilon\gamma}{R\sqrt{d}} \lesssim\frac{\varepsilon^2}{M_2R\sqrt{d}},
\end{equation}
and from \eqref{eq:var_noise}\pd{, in order to control the variance, the following condition for $\Delta$} \gea{holds} (i.e. $d\Delta^2/\gamma^2 \le M_2^2$):
\begin{equation}\label{eq:var_nosie}
\Delta\lesssim\frac{\gamma M_2}{\sqrt{d}}\lesssim \frac{\varepsilon}{\sqrt{d}}.
\end{equation}
\gea{one can} observe that the last inequality is less restrictive than \eqref{lon}. \pd{Thus,} \eqref{lon} describes the maximum possible level of noise. This level of noise corresponds to the part of the lower bound from \cite{risteski2016algorithms}: $\Delta \gtrsim \max\left\{\frac{\varepsilon^2}{ M_2R\sqrt{d}},\frac{\varepsilon}{d}\right\}$ obtained for $p=2$ (here and in other parts of the survey logarithmic factors in lower bounds \gea{are omitted}).

Note that in Sections~\ref{sec:direct} and ~\ref{sec:rand_smooth} only <<biased inequalities>> \gea{are considered} (see \eqref{eq:err}, \eqref{eq:unbias}) to obtain the maximum level of admissible noise. In Section~\ref{sec:non_smooth}, <<biased inequality>> \eqref{eq:unbias_non_smooth} \gea{is also considered} (that gives \eqref{lon}), but \gea{one also has} to consider inequality \eqref{eq:var_noise} (that gives \eqref{eq:var_nosie}). If $\varepsilon$ is small enough \eqref{lon} is more restrictive and \gea{one can} skip \eqref{eq:var_nosie}. \pd{In the same way}, the details at these parts in Sections~\ref{sec:direct},~\ref{sec:rand_smooth} \gea{are skipped}, assuming that $\varepsilon$ is sufficiently small. 

Note that the described above approach can be applied for smooth problems. For example, if $f$ has $L$-Lipschitz gradient in $2$-norm, then in Theorem~\ref{main_properties} \cite{gasnikov2016gradient,nesterov2017random} \gea{one should have}
 \begin{equation}\label{eq:hls}
 f(x) \leq f_{\gamma}(x) \leq f(x) + \frac{L\gamma^2}{2},
  \end{equation}
    $$\|\nabla f_{\gamma}(y) - \nabla f_{\gamma}(x) \|_2 \leq L \| y - x\|_2.$$ 
\pd{This allows one} to improve the number of iterations \pd{from} $d^{1/4}\varepsilon^{-1}$ \pd{to} $\varepsilon^{-1/2}$ (this is an expected result, see Sections~\ref{sec:direct},~\ref{sec:rand_smooth}), make the smoothing parameter large: $\gamma = \sqrt{\varepsilon/L}$, and, as a consequence, to make maximum admissible level of noise also larger (see \eqref{lon}): $$\Delta \le \frac{\varepsilon^{3/2}}{R\sqrt{dL}}.$$

If $f$ has higher level of smoothness\pd{, then a} more accurate analysis is required  since \eqref{eq:hls} does not have trivial higher-order generalizations. Moreover, to improve the oracle complexity \pd{from} $d\varepsilon^{-2}$ \pd{to} $d\varepsilon^{-1/2}$, it is required to use accelerated methods \textbf{A}($L,\sigma^2$) with <<strong growth condition>>  \cite{vaswani2019fast} rather than ordinary ones. \pd{This leads to a worse} iteration complexity. Thus, the improvement of the maximum admissible level of noise makes at least one of the other two criteria worse.

By using the restart technique \cite{nemirovsky1983problem,juditsky2014deterministic}, one can prove a counterpart of Theorem~\ref{FGM} for the case when $f$ is $\mu$-strongly convex w.r.t. the $p$-norm for some $p\in[1,2]$ and $\mu \ge \varepsilon/R^2$.
\begin{theorem}\label{SFGM}
Based on restarted batched Accelerated gradient method, the smoothing scheme applied to non-smooth and strongly convex problem \eqref{eq:problem}, provides a gradient-free method with
$\tilde{O}\left(\frac{d^{1/4}\sqrt{M M_2}}{\sqrt{\mu\varepsilon}}\right)$   successive iterations and   $\tilde{O}\left(\frac{\min \left\{ q, \ln d \right\} d^{\frac{2}{q}} M_2^2}{\mu\varepsilon}\right)$
zeroth-order oracle calls. 
\end{theorem}

Analogously, all other results in this paper (except the results from the end of Section~\ref{sec:non_smooth_one_point}, from Section~\ref{sec:online} and some results on maximum admissible level of noise) can be transferred for strongly convex problem \eqref{eq:problem}.

\subsection{$l_1$-ball smoothing scheme}\label{sec:non_smooth_l1}
In the core of the \dd{E}uclidean smoothing scheme lies the Stokes' formula \cite{akhavan2022gradient} (here $\nabla f(x)$ is subgradient!):
\begin{equation}\label{eq:Stoks}
    \int_{D} \nabla f(x) dx = \int_{\partial D} f(x)n(x)dS(x),
\end{equation}
where $\partial D$ is a boundary of $D$, $n(x)$ is a outward normal vector at point $x$ to $\partial D$, and $dS(x)$ denotes the surface measure. The Stokes' formula \gea{implies} $$\nabla f_{\gamma}(x)=  \mathbb{E}_{u} \left[\nabla_x f(x + \gamma u)\right] = \mathbb{E}_{e}\left[  d\frac{f(x+\gamma e)}{\gamma}e\right],$$ where $u$ \pd{is} \pd{random} uniformly distributed on $B_2^d(1)$ and  $e$ is \pd{random} uniformly distributed on $S_2^d(1)$. In \cite{gasnikov2016gradient} and independently in \cite{akhavan2022gradient}, it was mentioned that from the Stokes' formula \gea{one can} also obtain
$$\nabla \tilde{f}_{\gamma}(x)=  \mathbb{E}_{u} \left[\nabla_x f(x + \gamma \tilde{u})\right] = \mathbb{E}_{\tilde{e}}\left[  d^{3/2}\frac{f(x+\gamma \tilde{e})}{\gamma}n(\tilde{e})\right],$$
where $\tilde{u}$ \pd{is} \pd{random} uniformly distributed on $B_1^d(1)$,  $\tilde{e}$ is \pd{random} uniformly distributed on $S_1^d(1)$, and $n(\tilde{e})=\frac{1}{\sqrt{d}}\text{sign}(\tilde{e})$.

\pd{Thus,} \gea{one can} \pd{introduce a smoothed version of $f$}
$$
\tilde{f}_{\gamma}(x)=  \mathbb{E}_{\tilde{u}} f(x + \gamma \tilde{u}),
$$
where $\tilde{u} \sim RB^d_1(1)$, i.e., $\tilde{u}$ is random vector uniformly distributed on a unit ball in $1$-norm \pd{denoted as} $B^d_1(1)$. The following theorems are taken from \cite{gasnikov2016gradient,akhavan2022gradient,lobanov2022}.

\begin{theorem}[properties of $\tilde{f}_{\gamma}$]
\label{main_properties_l1} For all $x,y\in Q$
\begin{itemize}[noitemsep,topsep=0pt,leftmargin=10pt]
    \item $f(x) \leq \tilde{f}_{\gamma}(x) \leq f(x) + \frac{1}{\sqrt{d}}\gamma M_{2}$;
    \item $\tilde{f}_{\gamma}(x)$ has $L = 
    \dfrac{d M}{\gamma} $-Lipschitz gradient:
    $$\|\nabla \tilde{f}_{\gamma}(y) - \nabla \tilde{f}_{\gamma}(x) \|_{\avg{q}} \leq L \| y - x\|_{\avg{p}},$$
    where $|f(y) - f(x) | \leq M \|y - x \|_{\avg{p}}$.
\end{itemize}
\end{theorem}

\pd{Further, a stochastic gradient} \gea{is introduced as follows:}
\begin{equation}\label{sg_l1}
    \nabla \tilde{f}_{\gamma}(x,\tilde{e}) = d\frac{f(x+\gamma \tilde{e}) - f(x-\gamma \tilde{e})}{2\gamma} \text{sign}(\tilde{e}),
\end{equation}
where $\tilde{e} \sim RS^d_1(1)$ is a random vector uniformly distributed on a unit sphere in $1$-norm \pd{denoted as} $S^d_1(1)$. In comparison with \eqref{sg}, $\tilde{f}_{\gamma}(x,\tilde{e})$ requires less memory ($1$ float and $d$ bits). \pd{For simplicity} in this formula and in the next theorem,  \gea{it is assumed} that $\Delta = 0$.

    \begin{theorem}[properties of $\nabla \tilde{f}_{\gamma}(x,\tilde{e})$]
\label{main_properties2_l1} 
For all $x\in Q$
 \begin{itemize}[noitemsep,topsep=0pt,leftmargin=10pt] 
    \item $\nabla \tilde{f}_{\gamma}(x,\tilde{e})$ is an unbiased approximation for $\nabla \tilde{f}_{\gamma}(x)$:
    $$\mathbb{E}_{\tilde{e}} \left[\nabla \tilde{f}_{\gamma}(x, \tilde{e}) \right] = \nabla \tilde{f}_{\gamma}(x);$$
    
    \item  $\nabla \tilde{f}_{\gamma}(x,\tilde{e})$ has bounded variance (second moment):
   \begin{equation*}\label{eq:var_l1}
    \mathbb{E}_{\tilde{e}} \left[ \| \nabla \tilde{f}_{\gamma}(x, \tilde{e}) \|^2_q \right] \le 288 d^{\frac{2}{q}} M_2^2.
    \end{equation*}
   \end{itemize}
   \end{theorem}

\pd{Summarizing,} the results of Section~\ref{sec:non_smooth_l2} 
\avg{take} place for \dd{the} $l_1$-smoothing scheme, but with potentially better (for huge $d$ -- in practice this potential advantage \gea{is not observed} \cite{lobanov2022}) factor:
\avg{$2\sqrt{2} \min \left\{ q, \ln d \right\}\to 288$, $q\ge 2$}

\pd{Such randomization leads to the number of successive iterations bounded as}
$$O\left(\frac{d^{1/4}\sqrt{M M_2}R}{\varepsilon}\right), $$
\pd{and} the number of oracle calls \pd{bounded as} (may depend on $d$ only through the constants  \avg{$R = O\left(\|x^0 - x_*\|_p\sqrt{\min\left\{q,\ln d\right\}}\right)$ and $M_2$):}
$$O\left(\frac{ M_2^2 R^2}{\varepsilon^2}\right).$$
The level of maximum admissible noise coincides with \eqref{lon}: $\Delta \lesssim\frac{\varepsilon^2}{M_2R\sqrt{d}}$, where $R = \|x^0 - x_*\|_2$. Note that for $p=1$ the lower bound is not known. \gea{One} may expect that this bound could be further improved \cite{risteski2016algorithms} (i.e., maximum admissible level of noise could be large).

\subsection{Stochastic optimization problems. Two-point feedback}\label{sec:non_smooth_stoch}
 In this section, for certainty, \dd{the} $l_2$-smoothing scheme approach \gea{is considered}. The same can be done for any other approaches (e.g., \pd{for the algorithms corresponding to} Section~\ref{sec:rand_smooth}, see \cite{gasnikov2015accelerated,gorbunov2022accelerated}).

If the zeroth-order oracle returns an unbiased noisy stochastic  function value $f(x,\xi)$ ($\mathbb{E}_{\xi} f(x,\xi) = f(x)$), then with two-point oracle \cite{duchi2015optimal,gasnikov2017stochastic} \gea{one can} introduce the following counterpart of \eqref{sg}
\begin{equation}
\label{sg_stoch}
     \nabla f_{\gamma}(x,\xi,e) = d\frac{f(x+\gamma e,\xi) - f(x-\gamma e,\xi)}{2\gamma}e.
\end{equation}
\dd{In this case, a stochastic counterpart of Theorem~\ref{main_properties2} can be formulated by}
 the appropriate changes of $\nabla f_{\gamma}(x,e)$  to $\nabla f_{\gamma}(x,\xi,e)$, the expectation $\mathbb{E}_{e}$ to the expectation $\mathbb{E}_{e,\xi}$, and the redefinition of $M_2$ as a constant satisfying $\mathbb{E}_{\xi}\|\nabla_x f(x,\xi)\|_2^2\le M^2_2$ for all $x\in Q$.
\dd{Moreover, under this redefinition of  constant $M_2$ \gea{one can} also reformulate Theorem~\ref{FGM}.}
 
 Note that for high-probability deviation bounds in \dd{the} case of deterministic oracle \dd{returning the values of} $f(x)$, \gea{one can show} sub-\dd{G}aussian concentration \cite{dvinskikh2022gradient}. But there \pd{is} a lack of results for stochastic oracle \dd{returning the realizations of}  $f(x,\xi)$, see also \cite{dvinskikh2022gradient}. Fortunately, in  the approach \dd{in Section~\ref{sec:non_smooth_l2}} \gea{one can} take $\textbf{A}(L,\sigma^2)$ to be clipped batched Accelerated gradient method from \cite{gorbunov2021near} that guarantee almost sub-\dd{G}aussian concentration iff $\mathbb{E}_{\xi}\|\nabla_x f(x,\xi)\|_2^2\le M^2_2$. 

\subsection{Stochastic optimization problems. One-point feedback}\label{sec:non_smooth_one_point}
 In this section, for certainty, \dd{the} $l_2$-smoothing scheme approach \gea{is considered as a starting point}.

If the two-point feedback  \eqref{sg_stoch} is \dd{un}available, \dd{the} $l_2$-smoothing technique can utilize the one-point feedback by using the unbiased estimate \cite{nemirovsky1983problem,flaxman2005online,gasnikov2017stochastic}:
\begin{equation*}
    \nabla f_{\gamma}(x,\xi,e) = d\frac{f(x+\gamma e,\xi)}{\gamma}e
\end{equation*}
with \cite{gasnikov2017stochastic}
$$ \mathbb{E}_{\xi,e} \left[ \| \nabla f_{\gamma}(x, \xi, e) \|^2_q \right] \le   \begin{cases}
   \frac{(q-1)d^{1+2/q} G^2}{\gamma^2}, \;  q\in[2,2\ln d]\\
  \frac{4d(\ln d) G^2}{\gamma^2}, \; q \in (2\ln d,\infty),
 \end{cases}$$
where $\gamma$ is defined in \eqref{gamma} and it is assumed that  $\mathbb{E}_{\xi} \left[ | f(x, \xi) |^2\right] \le G^2$ for all 
$x\in Q$.
Thus, the \textit{Smoothing technique} can be generalized to the one-point feedback setup by replacing the RHS of \eqref{var} by the above estimate. This leads to the same iteration complexity $N$, but increases the oracle complexity $T$, and, consequently, the batch size  $B$ at each iteration.

In particular, for $p=2$ ($q=2$) $T\simeq d^2G^2M_2^2R^2/\varepsilon^4$ and $T\simeq d^2G^2M_2^2/(\mu\varepsilon^3)$ if $f$ is $\mu$-strongly convex in $2$-norm \cite{nemirovsky1983problem,flaxman2005online,gasnikov2017stochastic}. If $f$ has $L$-Lipschitz gradient in $2$-norm then \eqref{eq:hls} holds true and $\gamma = \sqrt{\varepsilon/L}$. Hence, $T\simeq d^2G^2LR^2/\varepsilon^3$ and $T\simeq d^2G^2L/(\mu\varepsilon^2)$ respectively \cite{agarwal2010optimal,gasnikov2017stochastic,akhavan2020exploiting}. The last bound corresponds to the lower bound \cite{jamieson2012query,shamir2013complexity}. For non-smooth convex problems the lower bound is similar $\gtrsim d^2/\varepsilon^2$ \cite{dani2007price,akhavan2020exploiting}, but it could be not tight. 
\dd{Here, for the first time, a situation where the strongly convex case and the convex case are not treated in the same way} \gea{is encountered}.

An important observation related \dd{to} one-point feedback is \pd{the} following: in some works (see e.g. \cite{gasnikov2014stochastic,gasnikov2015gradient,akhavan2020exploiting,akhavan2022gradient,akhavan2023gradient,lobanov2023black}) assumptions about adversarial noise $\delta$ \dd{imply} that $\delta$ does not depend on \pd{the used} randomization $e$ (and $\kappa$, see Section~\ref{sec:ZOsmooth}). In this case, $\E_e \left[\delta e\right] \equiv 0$ (see also \cite{granichin2003randomizirovannye}, where additional randomization reduce adversarial noise to stochastic noise) and in terms of convergence \dd{of function values} in \dd{the} expectation  \gea{one can} consider this noise $\delta$ to be stochastic unbiased. Moreover, if the \dd{the maximal} level of noise $\Delta \simeq G$ is small it make\dd{s} sense to split the complexities for two-point feedback and one-point stochastic unbiased feedback. \pd{This} was done in \cite{akhavan2020exploiting,akhavan2022gradient} without the acceleration techniques described in the  previous sections.

All the results of the survey can be improved in terms of the dependence o\dd{n} $\varepsilon$ by worsening the dependence on $d$. 
Namely, \pd{one} can use centre of gravity methods with iteration complexity $O(d\ln \varepsilon^{-1})$ \cite{nemirovsky1983problem} and estimate (sub-)gradient by using one-point feedback. For instance, \dd{the} Vayda's cutting plane method does not accumulate (sub-)gradient  error  \dd{over} iterations \cite{gladin2022vaidya}. For smooth problems with deterministic oracle, 
this is rather obvious. For example, in oracle  complexity bound from  Section~\ref{sec:direct} \gea{one can} replace  $d\varepsilon^{-1/2}$ by $d^2 \ln \varepsilon^{-1}$.
For non-smooth problems, such a
\dd{change}
$d\varepsilon^{-2}\to d^2 \ln \varepsilon^{-1}$ is not 
\dd{obvious} \cite{nemirovsky1983problem}, but still holds  under proper correction of the approach \cite{protasov1996algorithms}. For stochastic (one-point) feedback, the  best known oracle complexity result is  $O\left(d^{7.5}\varepsilon^{-2}\right)$ \cite{nemirovsky1983problem,agarwal2011stochastic,belloni2015escaping} (see also conjecture $O\left(d^3\varepsilon^{-2}\right)$ about the optimal bound in \cite{bubeck2017kernel}).

 If $f$ has $L$-Lipschitz gradient in \dd{the} $2$-norm, then \dd{the} Kernel-based approach from Section~\ref{sec:ZOsmooth} allows \pd{one} to improve oracle complexity \pd{to} 
  \al{$d^{2}\varepsilon^{-\frac{2(n+1)}{n}}$}, when $n$ is large
\cite{akhavan2020exploiting,novitskii2022improved,akhavan2023gradient}. Note that the lower bound (for $n \ge 2$) is \cite{akhavan2020exploiting}:
\[ \min\left(\dfrac{d^{1+\frac{1}{n}}}{\varepsilon^{\frac{2(n+1)}{n}}}, \dfrac{d^2}{\varepsilon^2}\right).\]
\dd{Hence, a significant gap still exists between lower bounds and current methods.}  

\subsection{Online optimization}\label{sec:online}
\dd{Nowadays, 
due to numerous applications in the  reinforcement learning, online optimization \cite{cesa2006prediction,hazan2016introduction,orabona2019modern} has become increasingly popular.}
\pd{Possibly,} the best known online problem formulation \pd{with zeroth-order oracle} is the multi-armed bandit problem  \cite{bubeck2012regret,slivkins2019introduction,lattimore2020bandit}, \pd{for which} the optimal bound $d\varepsilon^{-2}$ for the regret  is attained. For more general zeroth-order online convex optimization problems (on general convex sets, not necessarily \pd{on the} unit simplex as for multi-armed bandits), the lower bound is $d^2\varepsilon^{-2}$. This bound was obtained in the class of linear functions \cite{dani2007price}. The best known upper bound for general online convex optimization problem with one-point bandit feedback is $O\left(d^{9.5}\varepsilon^{-2}\right)$ \cite{bubeck2017kernel}.

For online setup, the iteration complexity typically coincides with the oracle complexity. Moreover, oracle complexity bounds for two-points online feedback are almost the same as for non-online case \cite{gasnikov2017stochastic} (up to logarithmic factors). Note that in \cite{gasnikov2017stochastic} the smoothness of $f$ was used to estimate the second moment of stochastic gradient. But, it can be easily replaced by \eqref{eq:var_noise} for non-smooth case. Also, note that the approach from Section~\ref{sec:non_smooth_l1} (see Theorem~\ref{main_properties2_l1}) can be applied to the online setting. This approach reduces $2\sqrt{2} \min \left\{ q, \ln d \right\}\to 288$, $q\ge 2$.

To conclude this section, \gea{it is important to note} that any advances in (online/stochastic) convex optimization usually involve advances in gradient-free \pd{optimization}. For example, based on \cite{zhang2022parameter} it seems \pd{to be} possible to build optimal adaptive gradient-free methods for convex stochastic online optimization problem with heavy-tails.

\section{Further research}\label{sec:further}
\pd{In the previous sections}, general convex optimization problem \eqref{eq:problem} \gea{are considered}. \dd{However,} many \dd{aforementioned} results (tricks)  can be naturally \pd{used for} saddle-point problems, sum-type problems, distributed optimization. 

\subsection{Saddle-point problems}\label{sec:SPP}
\gea{Consider} non-smooth convex-concave saddle-point problems
\begin{equation}\label{SPP}
    \min_{x\in Q_x\subseteq \mathbb{R}^{d_x}} \max_{y\in Q_y\subseteq \mathbb{R}^{d_y}} f(x,y).
\end{equation}
Gradient-free methods for convex-concave saddle-point problems were studied in \cite{beznosikov2020gradient,beznosikov2021one,gladin2021solving,sadiev2021zeroth,gasnikov22a,dvinskikh2022gradient,sadykov2023gradient}. 

Applying $l_2$-smoothing technique separately to $x$-variables and $y$-variables, \gea{one can} obtain almost the same results as for optimization problems with the only difference in Theorem~\ref{main_properties}: instead of
$$f(x) \leq f_{\gamma}(x) \leq f(x) + \gamma M_{2}$$ \gea{another inequality is considered:} $$f(x, y) - \gamma_y M_{2,y}\leq f_{\gamma}(x, y) \leq f(x, y) + \gamma_x M_{2,x}.$$
This leads to a clear counterpart of \eqref{gamma} for choosing $\gamma = \left(\gamma_x,\gamma_y\right)$, where $M_{2,x}$, $M_{2,y}$ -- corresponding Lipschitz constants in $2$-norm.

If \gea{one takes} as \textbf{A}($L,\sigma^2$) the batched Mirror-Prox or the batched Operator extrapolation method or the batched Extragradient method \cite{juditsky2011solving,kotsalis2020simple,gorbunov2021stochastic}, using  \eqref{L}, \eqref{var}, \gea{one can} obtain the following bounds
\begin{center}
$N(L,\varepsilon) = O\left(\frac{LR^2}{\varepsilon}\right)=O\left(\frac{\sqrt{d}M M_2R^2}{\varepsilon^2}\right),$ $ T(L,\sigma^2,\varepsilon) = O\left(\max\left\{N(L,\varepsilon),\frac{\sigma^2R^2}{\varepsilon^2}\right\}\right)$
$=O\left(\frac{\min \left\{ q, \ln d \right\} d^{\frac{2}{q}} M_2^2R^2}{\varepsilon^2}\right),$ 
\end{center}
where $d = \max\{d_x,d_y\}$,
$M_2=\max\left\{M_{2,x},M_{2,y}\right\}$, $R$ depends on the \pd{approximate optimality} criteria. For example, if $\varepsilon$ is expected accuracy in \pd{terms of the} fair duality gap \cite{juditsky2011solving}, then $R$ is a diameter in $p$-norm of $Q_x\otimes Q_y$ up to a $\sqrt{\ln d}$ factor, where $\otimes$ is the Cartesian product of two sets. The last equality assumes that $d \lesssim (M_2/M)^{2}$ when $p=1$. This result is also correct for stochastic saddle-point problems with proper redefinition of what is $M_2$, see Section~\ref{sec:non_smooth_stoch}.

Analogously, \gea{one can} take \textbf{A}($L,\sigma^2$) to be more specific batched algorithm that takes into account, for example, different constants of strong convexity and concavity \cite{gasnikov22a,metelev2022decentralized,li2022nesterov}.

Overall, all of the above results (except Sections~\ref{sec:l2} and \ref{sec:online}) can be  transferred to saddle-point problems. In particular, for the smooth case, a variety of optimal (stochastic) gradient-based algorithms for (stochastic) saddle\pd{-point} problems are collected in the survey \cite{beznosikov2022smooth}. All these algorithms could be extended to gradient-free oracle by using the tricks described above. But there are only few works (mentioned above) where such attempts have been made.

\subsection{Sum-type problems}\label{sec:sum-type}
For non-smooth sum-type problems, \dd{the} $l_2$-smoothing scheme (see Section~\ref{sec:non_smooth_l2}) can be applied, see \cite{gasnikov22a} for details. Moreover, all of the above results (except Section~\ref{sec:online}) can be  transferred to sum-type problems. Sometimes it requires non-trivial combinations of variance-reduced technique and coordinate-wise randomization technique \cite{hanzely2020variance}, see also Section~\ref{sec:coord}. But there is also a lack of papers in this field.

\subsection{Distributed optimization}\label{sec:distr}
Distributed optimization (on time-varying networks) is currently a burgeoning area of research, see, e.g., \cite{gorbunov2022recent,rogozin2022decentralized}. 

In distributed optimization, \gea{there is} one more \gea{criterion} (number of communications) that typically reduces to the number of successive iterations consider\dd{ed} before. So, for smooth problems, the describe\dd{d} above results (except  Sections~\ref{sec:coord},~\ref{sec:l2},~\ref{sec:online}) \dd{make it possible} to build optimal gradient-free distributed methods based on gradient ones. Unfortunately, \gea{there are no} references where this has been done.

For non-smooth stochastic problems, by using the Lan's sliding \cite{lan2020first}, in \cite{beznosikov2020derivative} an optimal algorithm was proposed with two-point feedback and in \cite{stepanov2021onepoint} the best-known algorithm  was proposed  with one-point feedback (see also \cite{akhavan2021distributed}). Note that \dd{the} $l_2$-smoothing technique from Section~\ref{sec:non_smooth_l2} gives non-optimal communication steps number, since the number of successive iterations $\sim d^{1/4}\varepsilon^{-1}$ while for optimal algorithms it should be $\sim \varepsilon^{-1}$ \cite{gasnikov22a}. That is why it required to use another technique \cite{beznosikov2020derivative,stepanov2021onepoint}.
\pd{For this direction \gea{one can} also observe a lack of papers.}

\section{Conclusions}

\gea{This survey discusses the state-of-the-art results and techniques in derivative-free convex optimization. The paper mainly focuses on the theoretical results about the convergence of the methods. In addition to the general convex minimization, the survey contains a discussion of the results for saddle-point problems, sum-type problems, and distributed optimization.}

\section{Cross-references}

\gea{\begin{itemize}
    \item Online Convex Optimization
    \item Decentralized Convex Optimization over Time-Varying Graphs
    \item Saddle Point Theory and Optimality Conditions
    \item Unified analysis of SGD-type methods
\end{itemize}}

\section{Acknowledgments}

The authors are grateful to Oleg Granichin, Katya Scheinberg, Alexander Tsybakov, Ivan Oseledets \avg{and Wouter Jongeneel}.

The work was supported by the Ministry of Science and Higher Education of the Russian Federation (Goszadaniye) 075-00337-20-03, project No. 0714-2020-0005.

The paper is dedicated to the memory of \textit{Yu. M.~Ermoliev (1936 -- 2022)} -- one of the founders of Soviet Stochastic Optimization school.

\bibliographystyle{abbrv}
\bibliography{references}

\end{document}